\documentclass{article}
\usepackage[english]{babel}
\usepackage{amsfonts, amsmath, amsthm, amssymb,amscd,indentfirst}
\usepackage{graphicx}
\usepackage{wrapfig}
\usepackage[rflt]{floatflt}
\usepackage{amsmath,amssymb,latexsym,indentfirst}
\usepackage{enumerate,pst-all}
\usepackage{makeidx}
\usepackage{color}
\usepackage{times}
\usepackage{palatino}

\newtheorem{theorem}{Theorem}[section]

\newtheorem{lemma}{Lemma}[section]

\newtheorem{remark}{Remark}[section]

\def\Cc{{\cal C}}

\def\ovn{\nabla}

\def\l{{\lambda}}

\def\l{{\lambda}}

\def\O{{\Omega}}

\def\Cc{{\cal C}}

\def\ovn{\nabla}

\def\Cc{\mathcal{C}}

\font\bbbld=msbm10 scaled\magstep1

\newcommand{\bfR}{\hbox{\bbbld R}}

\newcommand{\om}{\Omega}
\newcommand{\I}{\infty}

\begin{document}

\title{A theorem of do Carmo-Zhou type: oscillations and estimates for the first eigenvalue of the p-Laplacian}%

\author{Barnab\'e Pessoa Lima\thanks{Universidade Federal do Piau\'{\i},
 Departamento de Matematica, Campus Petronio Portela, Ininga, 64049-550 Teresina/PI, Brazil ({\tt barnabe@ufpi.br}). Partially supported by a PROCAD-CAPES fund} \and Newton Lu\'is Santos\thanks{Federal University of
Piau\'{\i},
 Departament of Mathematics, Campus Petronio Portela, Ininga, 64049-550 Teresina/PI, Brazil ({\tt newtonls@ufpi.br}). Partially supported by   PROCAD-CAPES fund}} \maketitle

\centerline{Dedicated to Prof. Keti Tenenblat}

\begin{abstract}
It is shown that the estimates obtained by Manfredo P. do Carmo and Detang Zhou, in their paper {\it Eigenvalue estimate on complete noncompact Riemannian manifolds and applications} ({\em Trans. Amer. Math. Soc.} 351 (1999),4, 1391-1401), for the first eigenvalue of the Laplace-Beltrami operator on open manifolds, via an oscillation theorem, can be naturally extended for the semi-elliptic singular operator operator, p-Laplace on manifolds, defined by $\Delta_p(u)=-div (\|\nabla u\|^{p-2} \nabla u)$.
\end{abstract}

\section{Introduction and statement of results}

Recently, an increasing movement towards understanding equations and technics modeling physical phenomena via tools from Differential Geometry is in process. Taking into account the geometrical data of the considered system has allowed new and interesting results (see \cite{CCFS}, \cite{Gr}, \cite{DGL}) relating analytical and geometrical aspects.

One of these problems is related to the divergent semi-elliptic singular operator, $p$-Laplacian, defined for $p\in (0,\infty)$ and $u$ smooth function defined on an open domain, $\om\subset \bfR^n$
\begin{equation}\label{p-lap} \Delta_p(u)=-div (\|\nabla u\|^{p-2} \nabla u)\end{equation}
More generally, when $\om$ is an open domain in a smooth Riemannian manifold, $(M,g)$, the p-Laplacian, $\Delta_p$, is defined as in (\ref{p-lap}), where $\nabla=\nabla^g$ is the Riemannian connection and $div=div_g$ the Riemannian divergence. This operator appears naturally from the variational problem associated to the energy functional
\begin{equation*}
E_p: W^{1,p}_0(\O)\to \bfR\qquad\mbox{given by}\qquad E_p(u)=\int_\O \|\nabla u\|^p\,\,\, d\O
 \end{equation*}
where $W^{1,p}_0(\O)$ denotes the Sobolev space given by the closure of ${\cal C}^{\infty}(\O)$ functions compactly supported in $\O$, for the norm
\begin{equation*}\|u\|_{1,p}^p=\int_\om |u|^p\,\,\, d\O+\int_\O \|\nabla u\|^p\,\,\, d\O\end{equation*}
Observe that, when $p=2$, $\triangle_2$ is just the Laplace-Beltrami operator.

The $p$-Laplacian has emerged from mathematical models in Physics (see \cite{Hu},  \cite{Di}, \cite{Li}). For each value of $p$, the operator $\Delta_p$ has appeared in a variety of physical fields, like fluid dynamics, in the study of flow through porous media, nonlinear elasticity, and glaciology.

Many efforts have been done towards understanding the general p-Laplace equation $\Delta_pu=f(x,u)$, on a domain of $\O \subset \bfR ^n$ with smooth boundary and Dirichlet condition, for several classes of functions $f$.

Let $(M,g)$ be a smooth Riemannian manifold and $\om\subset M$ a domain. For any $p\in (1,\infty)$, the $p$-Laplacian on $\om$, $\Delta_p(u)=-div (\|\nabla u\|^{p-2} \nabla u)$, appears naturally from the variational problem associated to the energy functional
\begin{equation*}
E_p: W^{1,p}_0(\O)\to \bfR\qquad\mbox{given by}\qquad E_p(u)=\int_\O \|\nabla u\|^p\,\,\, d\O
 \end{equation*}
  where $W^{1,p}_0(\O)$ denotes the Sobolev space given by the closure of ${\cal C}^{\infty}(\O)$ functions compactly supported in $\O$, for the norm
  $$\|u\|_{1,p}^p=\int_\om |u|^p\,\,\, d\O+\int_\O \|\nabla u\|^p\,\,\, d\O$$
Observe that, when $p=2$, $\triangle_2$ is just the Laplace-Beltrami operator.
We are interested in the nonlinear eigenvalue problem
\begin{equation}
\triangle_p u+\lambda |u|^{p-2}u=0 \label{eigenv problem}
\end{equation}
Since solutions for this problem, for arbitrary  $p\in(1,\infty)$ are only locally ${\cal C}^{1,\alpha}$ (exceptions for the case $p=2$), they must be described in the sense of distribution, that is, $u\in W^{1,p}_0(\O)$, not identically $0$ is an eigenfunction, associated to the eigenvalue $\lambda$, if
$$
\int_\O \|\nabla u\|^{p-2}g(\nabla u,\nabla \phi)\,\,\, d\O =\lambda \int_\O |u|^{p-2} u\phi \,\,\, d\O
$$
for any test function $\phi\in \Cc^{\infty}_0(\O)$.

It is remarkable that $\Delta_p$ exhibits some very interesting analogies with the Laplacian. Thus, a reasonable question, although highly nontrivial, to point out is:

\medskip

\noindent{\bf General question:} What known results, estimates and technics holding for the linear theory of the Laplace-Beltrami operator can be extended for $\Delta_p$?

\medskip

Some preliminary and interesting results have been obtained in this direction: e.g. Holopainen {\it et all} \cite{HMP}, Gol'dshtein and Troyanov \cite{GT}, Lima {\it et all} \cite{LMS}, Ly \cite{Ly}, and Matei \cite{Ma1}.

In this paper we extend results obtained by Manfredo P. do Carmo and  Detang Zhou (\cite{CZ}), on estimation of the first eigenvalues for the Laplacian on complement of domains, precisely, we extend theorems 2.1 and 3.1 of \cite{CZ}. The main tool used there is an oscilation theorem and the relationship between a Liouville equation an associated Ricatti equation. Such a theorem can be naturally extended to the setup related to the $p$-Laplacian:

Let $p>1$ and denote $\Phi(s)=|s|^{p-2}s$. Below, $p$ and $q$ will represent conjugate exponents, that is $\dfrac1p+\dfrac1q=1$

Consider the Liouville type equation
\begin{equation}\label{assoc eq}
(v(t)\Phi(x'(t)))'+\l v(t)\Phi(x(t))=0
\end{equation}
This equation is intimately related to the following Ricatti equation:
\begin{equation}\label{Ricatti eq}
y'=\l v(t)+(p-1)v(t)^{1-q}y^q
\end{equation}
Indeed, if $x(t)$ is a positive solution for (\ref{assoc eq}), one set
\begin{equation}\label{Ricatti subst}
y(t):=-\dfrac{v(t)\Phi(x'(t))}{\Phi(x(t))}
\end{equation}
then, the first derivative of (\ref{Ricatti subst}) gives:
\begin{eqnarray*}
y'(t)&=&-\dfrac{(v(t)\Phi(x'(t)))'}{\Phi(x(t))}+\dfrac{v(t)\Phi(x'(t))\Phi'(x(t))}{\Phi^2(x(t))}\\
&=& \l v(t)+\dfrac{v(t)\Phi(x'(t))(p-1)x(t)^{p-2}x'(t)}{\Phi^2(x(t))}\\
&=& \l v(t)+(p-1)v(t)\dfrac{|x'(t)|^p}{|x(t)|^p}\\
&=& \l v(t)+(p-1)v(t)\dfrac{(|x'(t)|^{p-1})^{\frac{p}{p-1}}}{(|x(t)|^{p-1})^{\frac{p}{p-1}}}\\
&=&\l v(t)+(p-1)\dfrac{y^q}{v(t)^{q-1}}
\end{eqnarray*}
where, in the second equality we used equation (\ref{assoc eq}).

We say that the equation (\ref{assoc eq}) is oscillatory if a solution $x(t)$ for (\ref{assoc eq}) has zeroes in $[T,\I)$ for any $T\ge T_0$ (that is $x(t)$ has zeroes for arbitrarily large $t$). It is a classical result (see for instance \cite{Do}) that if one solution has arbitrarily large zeroes then, any other solution will also possess arbitrarily large zeroes, that is to be oscillatory is a characteristic of the equation and not to the particular solution.

Our main oscillation result is:

\begin{theorem}[compare theorem 2.1 \cite{CZ}]\label{thm1} Let $v(t)$ be a positive continuous function on $[T_0,+\infty)$, and define
\begin{equation*}
V(t):=\int^t_{T_0} v(s)\,\, ds, \qquad \mbox{and}\qquad \theta_V:=\lim_{t\to\I} \dfrac{\log V(t)}t
\end{equation*}
If the limit
\begin{equation*}\lim_{t\to\I}V(t)= \int^\infty_{T_0} v(s)\,\, ds =+\infty\end{equation*}
Then, equation (\ref{assoc eq}) is oscillatory provided that
\begin{description}
\item[(a)] $\lambda>0\quad $ and $\quad \theta_V=0$
\item or else
\item[(b)] $\lambda>\dfrac{c^p}{p^p}\quad $ and $\quad\theta_V\leq c$
\end{description}
\end{theorem}
\noindent As a particular case of item {\bf (a)}, if $\lambda>0$ and $V(t)\leq At^c$, for certain positive constants, $A$, $c$, then
$$\theta_V=\lim_{t\to\I} \dfrac{\log V(t)}t\leq \lim_{t\to\I} \dfrac{\log At^c}t=\lim_{t\to\I} \left(\dfrac{\log A}t+\dfrac{c\log t}t\right)=0$$
And, as a particular case of item {\bf (b)}, if $\lambda>\dfrac{c^p}{p^p}$ and $V(t)\leq Ae^{ct}$, then
$$\theta_V=\lim_{t\to\I} \dfrac{\log V(t)}t\leq \lim_{t\to\I} \dfrac{\log Ae^{ct}}t=\lim_{t\to\I} \left(\dfrac{\log A}t+\dfrac{ct}t\right)=c$$

\smallskip

Let $(M^n,g)$ be an open (that is complete, noncompact) Riemannian manifold and fix some base point $p\in M$. Denote by $v(r):= Area (\partial B_r(p))$ the area of the geodesic sphere of radius $r$ centered at $p$. Put
\begin{equation*}
V(r):=Vol(B_r(p))=\int_0^rv(r)\,\, dr \qquad\mbox{and set}\qquad \theta(M):=\lim_{t\to\I} \dfrac{\log V(t)}t
\end{equation*}
Notice that the number $\theta(M)$ does not depend on the base point, $p$, so it is an invariant of the manifold that captures the growth behavior $M$ off compact sets, since for any fixed $R>0$  holds
\begin{equation*}
\theta(M)=\lim_{t\to\I} \dfrac{\log (V(R)+A(R,t))}t=\lim_{t\to\I} \dfrac{\log A(R,t)}t
\end{equation*}
where $A(R,t)=V(t)-V(R)$. If $M$ is a manifold with $Vol (M)=+\I$ then, for any $T_0>0$, one has $\int_{T_0}^{+\I}v(r)\,\, dr=+\I$. Our main theorem is the following:

\begin{theorem}[compare theorem 3.1 of \cite{CZ}]\label{thm2} Let $(M,g)$ be an open manifold with infinite volume, $Vol(M^n)=+\I$. If $\O\subset M$ is an arbitrary compact set, denote by  $\lambda_{1,p}(M-\O)$ be the first eigenvalue for the $p$-Laplacian on $M-\O$
\begin{description}
\item[(a)] If $M$ has subexponential growth, that is $\theta(M)=0$, then $\lambda_{1,p}(M-\O) =0$
\item[(b)] If $\theta(M)\leq \beta$, for some $\beta>0$, then $\lambda_{1,p}(M-\O) \leq \dfrac{\beta^p}{p^p}$
\end{description}
\end{theorem}

\begin{remark}
As particular cases for this theorem we have:
\begin{enumerate}
\item If $Vol (B_r(p))\leq cr^\beta$ for some positive constants $c,\beta>0$ and $r\ge r_0$, then $\lambda_{1,p}(M-\O) =0$.
\item Let $(M^n,o,g,k)$ denote a pointed Riemannian manifold, with a base point, $o\in M$, and $k:[0,\I)\to \bf R$ a positive, nonincreasing function, with
    $$b_0(k):=\int_0^\I sk(s)\,\, ds<\I$$
    such that the radial Ricci curvatures (that is, the Ricci curvatures along radial direction, $\partial_r$, tangent to minimal unitary geodesics departing from the base point) satisfy the bound
    $${\rm Ric}_{_M}(\partial_r|_p)\ge -(n-1)k(dist(o,p)),\qquad \mbox{for all $p\in M$}$$
    we say that such a manifold has asymptotically nonnegative radial Ricci curvature, thus since for such a manifold one has $V(r)\leq e^{b_0(k)}r^n$ (see \cite{Zh}), it follows, by the previous case that $\lambda_{1,p}(M-\O) =0$
\item If $\,\, Vol (B_r(p))\leq ce^{\beta t}$, then $\lambda_{1,p}(M-\O) \leq \dfrac{\beta^p}{p^p}$
 \end{enumerate}
\end{remark}

\section{Proving the main results}

\noindent{\bf Proof of the theorem \ref{thm1}}
\smallskip

\noindent We shall assume that equation (\ref{assoc eq}) is nonoscillatory and this will lead us to a contradiction.

Suppose that equation (\ref{assoc eq}) is nonoscillatory and let $x(t)$ be any solution. Then since $\pm x(t)$ is solution for (\ref{assoc eq}) we can assume the existence of a certain $T\ge T_0$ such that $x(t)>0$ for all $t\in [T,\I)$ and let $y(t)$ the associated solution for the Ricatti equation (\ref{Ricatti eq}) given by the Ricatti substitution defined in (\ref{Ricatti subst}). Thus it is immediate that if $y(T)> 0$, $y'(t)>0$ for $t\in [T,+\I)$. As a consequence, $y$ is increasing in $[T,+\I)$. We have two cases to consider:
\begin{equation}\label{4}
\mbox{\bf (I)} \quad \int_T^{+\I} \frac{dr}{v(r)^{q-1}} <+\I \qquad\mbox{and}\qquad \mbox{\bf (II)}\quad \int_T^{+\I} \frac{dr}{v(r)^{q-1}} =+\I
\end{equation}

\noindent{\bf Case I.}

\noindent From (\ref{assoc eq}) we have $y'(t)\ge \lambda v(t)$ thus $y(t)\ge y(T)+\lambda (V(t)-V(T))$ is positive, for sufficiently large $t$, since $\lim_{t\to \I} V(t)=\I$. Consequently one can assume that $y(t)>0$ for all $t\ge T$.

\noindent Applying Young inequality, $AB\leq \dfrac{A^p}p+\dfrac{B^q}q$, if $A,B\ge 0$, for the couple $A=(p\l)^{1/p}v(t)^{q-1}$ and $B=(p-1)^{1/q}q^{1/q}y(t)$ implies that
\begin{equation*}
[(p\l)^{1/p}v(t)^{q-1}][(p-1)^{1/q}q^{1/q}y(t)]\leq \dfrac{(p\l)v(t)^{p(q-1)}}p+\dfrac{(p-1)qy(t)^q}q
\end{equation*}
Now, since $(p\l)^{1/p}(p-1)^{1/q}q^{1/q}=p\l^{1/p}$ we have
\begin{equation}\label{ineq1}
p\l^{1/p}v(t)^{q-1}y(t)\leq \l v(t)^q+(p-1)y(t)^q
\end{equation}
dividing inequality (\ref{ineq1}) by $v(t)^{q-1}$:
\begin{equation}\label{ineq2}
p\l^{1/p}y(t)\leq \l v(t)+(p-1)v(t)^{1-q}y(t)^q
\end{equation}
Comparing inequality (\ref{ineq2}) with Ricatti equation (\ref{Ricatti eq}) we get
\begin{equation}
y'(t)\ge p\l^{1/p}y(t)
\end{equation}
which after integration in $[T,t]$ gives
\begin{equation}\label{ineq sol}
y(t)\ge y(T)\exp (p\l^{1/p}(t-T))
\end{equation}
On the other hand Ricatti equation (\ref{Ricatti eq}) also implies the inequality $$y'\ge (p-1)v(t)^{1-q}y^q$$ or, else
\begin{equation}
\frac{y'(t)}{y(t)^q}\ge \frac{(p-1)}{v(t)^{q-1}}
\end{equation}
But, since $$\frac1{1-q}\frac{d}{dt}y(t)^{1-q}=\frac{y'(t)}{y(t)^q}$$ and $q>1$, it follows that
\begin{equation}
(1-q)\frac{y'(t)}{y(t)^q}=\frac{d}{dt}\left( y(t)^{1-q}\right) \leq (1-q)(p-1)v(t)^{1-q}=-v(t)^{1-q}
\end{equation}
Now, integrate this inequality over the interval $[t,t+1]$ to get
\begin{equation}\label{1}
y(t)^{(1-q)}\leq y(t-1)^{(1-q)}-\int_{t-1}^t \frac{ds}{v(s)^{q-1}}
\end{equation}
On the other hand, from (\ref{ineq sol}) we see that
\begin{equation}\label{2}
\frac1{y(t)^{q-1}}\leq \frac1{y(T)\exp (q\l^{1/p}(t-T))}
\end{equation}

H\"older inequality implies that
\begin{eqnarray*}
1&=&\int_{t-1}^t \frac{v(s)^{1/p}}{v(s)^{1/p}}\,\, ds\leq \left(\int_{t-1}^t (v(s)^{1/p})^p\,\, ds\right)^{1/p}\left(\int_{t-1}^t \frac1{v(s)^{q/p}}\,\, ds\right)^{1/q}\\
&=&\left(\int_{t-1}^t v(s)\,\, ds\right)^{1/p}\left(\int_{t-1}^t \frac1{v(s)^{q-1}}\,\, ds\right)^{1/q}\leq V(t)^{1/p}\left(\int_{t-1}^t \frac1{v(s)^{q-1}}\,\, ds\right)^{1/q}
\end{eqnarray*}
consequently
\begin{equation}\label{3}
\int_{t-1}^t \frac1{v(s)^{q-1}}\,\, ds\ge \frac1{\left(\int_{t-1}^t v(s)\,\, ds\right)^{q-1}} \ge \frac1{V(t)^{q-1}}
\end{equation}
If {\bf (a)} and the condition {\bf (I)} are satisfied, from (\ref{1}), (\ref{2}) and (\ref{3}) it follows that, for $t$ large enough:

\begin{eqnarray*}
0<\dfrac1{y(t)^{(q-1)}}&\leq &\frac1{y(T)\exp (q\l^{1/p}(t-T-1))}-\frac1{V(t)^{q-1}}\\
&\leq &\left(\frac{V(t)^{q-1}}{y(T)\exp (q\l^{1/p}(t-T-1))}-1\right)\frac1{V(t)^{q-1}}\\
&\leq &\left(\frac{\Lambda}{\exp t\left\{q\l^{1/p}-(q-1)t^{-1}\log V(t)\right\}}-1\right)\frac1{V(t)^{q-1}}<0
\end{eqnarray*}
a contradiction, where $\Lambda=\Lambda(q,\lambda,T):=\dfrac{\exp[q\l^{1/p}(T+1)]}{y(T)}$.

On the other hand, if {\bf (b)} and the condition {\bf (I)} are satisfied, a similar reasoning gives, for large enough $t$:
\begin{eqnarray*}
0<\dfrac1{y(t)^{(q-1)}}&\leq &\frac1{y(T)\exp (q\l^{1/p}(t-T-1))}-\frac1{c^{q-1}e^{(q-1)\alpha t}}\\
&\leq &\left(\frac{c^{q-1}\Lambda}{\exp \left\{(q\l^{1/p}-(q-1)\alpha )t\right\}}-1\right)\frac1{c^{q-1}e^{(q-1)\alpha t}}<0
\end{eqnarray*}
again a contradiction, since due to {\bf (b)}, $q\l^{1/p}-(q-1)\alpha=q(\l^{1/p}-\alpha/p)>0$.

\medskip

\noindent{\bf Case II.}
\smallskip

This case is a straightforward consequence of the following Leighton-Wintner type oscilation criteria, whose proof can be found in the interesting book of D\v{o}sl\'y and \v{R}eh\'ak, {\it Half linear differential equation} \cite{Dy}. As before, we let $p,q>1$ be conjugate exponents, and for $y\in \bfR$, denote $\Phi(y)=|y|^{p-2}y$:
\begin{lemma}[theorems 1.1.1 and 1.2.9 of \cite{Dy}]
Let  $r,c:I\to\bfR$ be continuous functions on some interval, $I\subset \bfR$, with $r(t)>0$ on $I$ and assume that $r\Phi (x')$ is of class ${\cal C}^1$. Fix $t_0\in I$, and numbers $A,B\in \bfR$. The half-linear second order equation
\begin{equation}\label{general Liouville}
(r(t)\Phi(x'))'+c(t)\Phi(x)=0
\end{equation}
has a unique solution on $I$ satisfying the initial conditions $x(t_0)=A$ and $x'(t_0)=B$. Moreover if $I$ has no upper bound, then equation (\ref{general Liouville}) is oscillatory provided
\begin{equation*}
\int^{\I}r^{1-q}(t)\,\, dt=+\I \qquad \mbox{and} \qquad \int^{\I}c(t)\,\, dt=\lim_{b\to\I}\int^{b}c(t)\,\, dt=+\I
\end{equation*}
\end{lemma}

\noindent In our case $r(t)=v(t)$ and $c(t)=\lambda v(t)$. With this lemma at hand, taking into account the hypothesis of theorem \ref{thm1}, and the condition assumed in {\bf Case II} (\ref{4}), we conclude the proof of theorem \ref{thm1}. $\square$

\bigskip

\noindent{\bf Proof of the theorem \ref{thm2}}
\smallskip

\noindent To prove this theorem we follow the steps of do Carmo and Zhou, in \cite{CZ}. Let
\begin{equation} \label{hip}
v(r)=Area(\partial B_r(p))\qquad\mbox{and denote}\quad B(r)=B_r(p)
 \end{equation}
If the volume $V(M)=+\I$, then $\int_T^{\I} v(t)dt=+\I$ for any constant $T>0$. Since $\Omega$ is a compact subset, we can find a constant $T_0$ such that $\Omega\subset B(T_0)$

\smallskip

If the condition {\bf (a)} is satisfied then from theorem \ref{thm1}, for any $\lambda$ there exists a nontrivial oscillatory solution
Thus there exists numbers $T_1^\lambda$ and $T_2^\lambda$ in $[T_0,+\I)$ with $T_1^\lambda<T_2^\lambda$ and $x_\l(T_1^\l)=x_\l(T_2^\l)=0$, and $x_\l(t)\neq 0$ for any $t\in (T_1^\l,T_2^\l)$

Set $\phi_\l(x)=x_\l (r(x))$, where $r(x):=dist_M(p,x)$ is the distance function from some $p\in M$. Denote by $\O_\l$ the annulus:
\begin{equation*}
\O_\l =B(T_2^\l)-B(T_1^\l) \subset M-\O
\end{equation*}
 Thus
$$
\lambda_{1,p}(M-\O)\leq \lambda_{1,p}(\O_\l)\leq \dfrac{\int_{\O_\l} |\ovn \phi_\l|^p\,\, dM}{\int_{\O_\l} |\phi_\l|^p\,\, dM}
$$
since $\ovn \phi_\l=x_\l '(r)\ovn r$ it follows that $|\ovn \phi_\l|^p=|x_\l '(r)|^p=|x_\l '(r)|^{p-2}x_\l '(r)x_\l '(r)$. Thus $|\ovn \phi_\l|^p=\Phi(x_\l '(r))x_\l '(r)$, consequently
\begin{equation}
\lambda_{1,p}(M- \O)\leq \dfrac{\int_{T_\l^1}^{T_\l^2}[v(r) \Phi(x_\l'(r))]x_\l'(r)\,\, dr}{\int_{T_\l^1}^{T_\l^2} |\phi_\l(r)|^pv(r)\,\, dr}=\lambda
\end{equation}
since $\l$ is an arbitrary positive constant, one has $\l_{1,p}(M-\O)=0$

\smallskip

 If the condition {\bf (b)} is satisfied then from theorem \ref{thm1} for any $\lambda>\frac{\beta^p}{p^p}$ there exists a nontrivial oscillatory solution of (\ref{assoc eq}) on $[T_0,+\I)$, with $v(t)$ as in (\ref{hip}). Thus there exists $T_1^\lambda$ and $T_2^\lambda$ in $[T_0,+\I)$ with $T_1^\lambda<T_2^\lambda$ and $x_\l(T_1^\l)=x_\l(T_2^\l)=0$, and $x_\l(t)\neq 0$ for any $t\in (T_1^\l,T_2^\l)$

As before we get:
\begin{eqnarray*}
\lambda_{1,p}(M-\O)&\leq \lambda_{1,p}(\O_\l)\leq \dfrac{\int_{\O_\l} |\ovn \phi_\l|^p\,\, dM}{\int_{\O_\l} |\phi_\l|^p\,\, dM} \leq \dfrac{\int_{T_\l^1}^{T_\l^2}[v(r) \Phi(x_\l'(r))]x_\l'(r)\,\, dr}{\int_{T_\l^1}^{T_\l^2} |\phi_\l(r)|^pv(r)\,\, dr}
\end{eqnarray*}
Since $\lambda$ is arbitrary positive constant larger than $\dfrac{\beta^p}{p^p}$, so $\lambda_{1,p}(M\setminus-\O)\leq \dfrac{\beta^p}{p^p}$. This proves part {\bf (b)} of our theorem. $\square$

\section{Applications}

In \cite{LMS} it is presented the concept of essential $p$-first eigenvalue associated to the $p$-Laplacian, and defined as
 \begin{equation}
\lambda_p^{ess}(M):=\lim_{j\to\infty} \lambda_{1,p}(M-K_j)
\end{equation}
where $K_1\subset K_2\subset \ldots$ is any exhaustion of $M$ through compact subsets - this limit being independent on the exhaustion (see \cite{Br}, \cite{Do} or \cite{DL}). This definition is a natural generalization of the spectral invariant associated to the Laplace-Beltrami operator, the essential spectrum, $\lambda_1^{ess}(M)$ which is the greatest lower bound of the essential spectrum of $M$, and admits the characterization $\lambda_1^{ess}(M)=\lim_{j\to\infty}\lambda_1(M- K_j)$, as above, where $\lambda_1(M- K_j)$ is the first eigenvalue for $\Delta$, on $M- K_j$.

\noindent Since item {\bf (b)} of theorem \ref{thm2}, implies that $\lambda_{1,p}(M-\O)\leq \dfrac{\theta(M)^p}{p^p}$ for any compact set $\O\subset M$ we re-obtain theorem {\bf 1.5} of \cite{LMS} (there obtained via a Brooks type technic, \cite{Br})

\begin{theorem}
If the volume of $M$ is infinity, then $\lambda_p^{ess}(M)\leq \frac{\theta(M)^p}{p^p}$.
\end{theorem}

\end{document}